\documentclass[12pt,a4paper]{article}
\usepackage{amsmath,amssymb,amsfonts}
\usepackage[latin1]{inputenc}

\def\Bbb{\mathbb}

\title{\bf The mean value of the digits of $1/p$}

\author{Kurt Girstmair}

\date{}

\makeatletter
\let\@@maketitle=\maketitle
\def\maketitle{\def\thispagestyle##1{\relax}\@@maketitle}
\makeatother
\newtheorem{theorem}{Theorem}

\def\BE{\begin{equation}}
\def\EE{\end{equation}}
\def\BD{\begin{displaymath}}
\def\ED{\end{displaymath}}
\def\BA{\begin{array}}
\def\EA{\end{array}}
\def\BEA{\begin{eqnarray*}}
\def\EEA{\end{eqnarray*}}
\def\BI{\bibitem}

\def\Q{\Bbb Q}

\def\phi{\varphi}

\def\MB{\mbox}
\def\LD{\ldots}

\def\SP#1{\langle #1 \rangle}

\def\DIV{\,|\,}
\def\NDIV{\, \nmid \,}

\def\MN{\medskip\noindent}

\def\BCHI{B_{\chi}}
\def\BCHIJ{B_{\chi^j}}

\newcommand{\btop}[2]{\genfrac{}{}{0pt}{1}{#1}{#2}}

\begin{document}
\maketitle
\normalsize

\begin{abstract}

Let $p\ge 3$ be a prime and $b\ge 2$ an integer such that $p$ does not divide $b$. Then $1/p$ has a periodic digit expansion with respect to the basis $b$. The length $l$ of the period
is the (multiplicative) order of $b$ mod $p$. If $l$ is even, then the mean value of the digits of a period is just $(b-1)l/2$. The case of an odd length  $l$ is more interesting.
If $l=(p-1)/2^m$ is odd, the mean value of the digits of a period was given previously.
This mean value involves generalized Bernoulli numbers. However, it is not clear how this result can be generalized to an arbitrary odd length $l$.
In the present note we settle this case.
\end{abstract}

\section*{1. Introduction and result}

Let $p\ge 3$ be a prime, $b\ge 2$ an integer such that $p\NDIV b$.
Then
\BE
\label{1.2}
  \frac 1p=\sum_{j=1}^{\infty}c_jb^{-j},
\EE
where the numbers $c_j\in\{0,1,\LD,b-1\}$ are the {\em digits} of $1/p$ with respect to the basis $b$.

For an integer $k$ let $(k)_p$ denote the integer $j\in\{0,\LD,p-1\}$ that satisfies $j\equiv k$ mod $p$. Then the digit $c_j$ of (\ref{1.2}) is given by
\BE
\label{1.4}
  c_j=\frac{b(b^{j-1})p-(b^j)_p}p, \: j\ge 1;
\EE
see \cite{Gi1}. This identity shows that
the sequence of the digits is periodic and that $(c_1,\LD,c_l)$ is a period, $l$ being the (multiplicative) order of $b$ mod $p$. The mean value of the digits $c_1,\LD,c_l$ is
$(c_1+\LD+c_l)/l$. Hence computing this mean value is the same as computing the sum
\BD
 S=\sum_{j=1}^lc_j.
\ED
This is easy if $l$ is even. Indeed, since $l$ is the order
of $b$ mod p, we have
\BD
  b^{l/2}\equiv -1 \MB{ mod } p.
\ED
Because $(b^{l/2+j})_p=p-(b^j)_p$, the identity (\ref{1.4}) shows $c_{{l/2}+j}=b-1-c_j$, and, accordingly
\BD
  S=\sum_{j=1}^{l/2}(b-1)=\frac{b-1}2\cdot l.
\ED
The first case with an odd length $l$ was settled in \cite{Gi1}, namely, $l=(p-1)/2$.
Indeed, in this case we have, for $p>3$,
\BE
\label{1.8}
 S=\frac{b-1}2\cdot l-\frac{b-1}2\cdot h,
\EE
where $h$ is the class number of $\Q(\sqrt{-p})$.
In the paper \cite{Gi2} this result is generalized to the case $l=(p-1)/2^m$, $l$ odd. The respective formula involves generalized Bernoulli numbers.
However, the method applied in the said paper does not produce a formula in the case of an arbitrary odd length $l$.

In the present paper we consider just this case. Therefore, let $l$ be an odd divisor of $p-1$. We write $l=(p-1)/q$, where $q$ divides $p-1$. Since $l$ is odd, $q$  is even, and we have
$2\NDIV (p-1)/q$, i.e., $p\equiv q+1$ mod $2q$.

Let $\chi$ be a Dirichlet character mod $p$ of order $q$. Then $\chi$ is odd. Indeed, let $g$ be a primitive root mod $p$.
Then $\chi(-1)=1$ is the same as $\chi(g^{(p-1)/2})=\chi(g)^{(p-1)/2}=1$. Since the order of $\chi(g)$ is $q$, we obtain $q\DIV (p-1)/2$ and $2\DIV (p-1)/q=l$.

Let
\BD
  \BCHI=\frac 1p\sum_{k=1}^{p-1}k\chi(k)
\ED
be the generalized Bernoulli number (of order $1$) belonging to $\chi$. In the case $l=(p-1)/2$, we have $\BCHI=-h$, provided that $p>3$; see (\ref{1.8}).

\begin{theorem} 
\label{t1}

Let the odd number $l$ be the order of $b$ mod $p$, $l=(p-1)/q$. Let $\chi$ be a Dirichlet character mod p of order $q$.
Then the sum $S$ of the digits $c_1,\LD,c_l$ of $1/p$
takes the value
\BD
   S=\frac{b-1}2\cdot l +\frac{b-1}2E,
\ED
where $E$ is the mean value of the Bernoulli numbers belonging to the odd characters in the group $\SP{\chi}$,
i.e.,
\BE
\label{1.12}
  E=\frac 2q\sum_{\btop{j=1}{j\:{ odd}}}^{q-1}\BCHIJ.
\EE

\end{theorem} 

\MN
{\em Remark.} Upper bounds for generalized Bernoulli numbers are known, such as
$|B_{\psi}|\le \sqrt p\log p$ for all odd characters $\psi$ mod $p$; see \cite[Lemma 4.8, Th. 4.9]{Wa} and \cite[Lemma 8.5]{Na}.
Thereby, the identity (\ref{1.12}) gives
\BD
|E|\le \sqrt p\log p.
\ED
Effective versions of the asymptotic estimate
\BD
   |B_{\psi}|\ll \sqrt p\log \log p,
\ED
which holds under the Generalized Riemann Hypothesis,
are known; see \cite{La}. We will also show $E\ne 0$.

\section*{2. Proof of Theorem \ref{t1}}

Observe that the numbers $(b^j)_p$, $j=1,\LD,L$, just run through the set of integers $k\in\{1,\LD,p-1\}$ with $\chi(k)=1$.
Therefore, the identity  (\ref{1.4}) yields
\BD
  S=\frac bp\sum_{\btop{k=1}{\chi(k)=1}}^{p-1}k-\frac 1p\sum_{\btop{k=1}{\chi(k)=1}}^{p-1}(bk)_p.
\ED
Since the numbers $(bk)_p$ of the second sum run through the same set, we obtain
\BE
\label{2.2}
 S=\frac{b-1}p\sum_{\btop{k=1}{\chi(k)=1}}^{p-1}k.
\EE
Now we decompose the right-hand side of (\ref{2.2}) into a sum with an explicit value and an alternating sum.
We have $\chi(k)=\pm 1$ if, and only if $\chi^2(k)=1$. However, $\chi^2$ is an even character, hence
\BD
  \chi^2(k)=\chi^2(p-k), \: k=1,\LD,p-1.
\ED
This gives
\BD
  \sum_{\btop{k=1}{\chi^2(k)=1}}^{p-1}k=\sum_{\btop{k=1}{\chi^2(k)=1}}^{(p-1)/2}(k+p-k)=p\cdot|\{k;1\le k\le(p-1)/2,\chi^2(k)=1\}|.
\ED
But the intervals $(0,p/2)$ and $(p/2,p)$ contain the same number of integers $k$ with $\chi^2(k)=1$. Since
$|\{k;1\le k\le p-1,\chi^2(k)=1\}|=2(p-1)/q$, we obtain
\BE
\label{2.4}
  \sum_{\btop{k=1}{\chi^2(k)=1}}^{p-1}k=\frac{p(p-1)}q=p\cdot l.
\EE
We use the obvious identity
\BD
  \sum_{\btop{k=1}{\chi(k)=1}}^{p-1}k=\frac 12\left(\sum_{\btop{k=1}{\chi^2(k)=1}}^{p-1}k+
  \sum_{\btop{k=1}{\chi^2(k)=1}}^{p-1}k\chi(k)\right)
\ED
in order to obtain, from (\ref{2.2}) and (\ref{2.4}),
\BE
\label{2.6}
 S=\frac{b-1}2\cdot l+\frac{b-1}2\cdot \frac 1p\sum_{\btop{k=1}{\chi^2(k)=1}}^{p-1}k\chi(k).
\EE
Observe that this identity is more fruitful than (\ref{2.2}). Indeed, in the case $l=(p-1)/2$, it is the same as (\ref{1.8}).
Moreover, there is some hope to connect the sum in (\ref{2.6}) with Bernoulli numbers. We are going to show
\BE
\label{2.8}
  \frac 1p\sum_{\btop{k=1}{\chi^2(k)=1}}^{p-1}k\chi(k)=E,
\EE
with $E$ as in (\ref{1.12}). We have
\BD
   \sum_{\btop{j=1}{j\:{ odd}}}^{q-1}\BCHIJ=\frac 1p\sum_{\btop{k=1}{\chi^2(k)=1}}^{p-1}k\chi(k)\cdot\frac q2+
   \frac 1p\sum_{\btop{k=1}{\chi^2(k)\ne 1}}^{p-1}k\sum_{\btop{j=1}{j\:odd}}^{q-1}\chi(k)^j.
\ED
But the sum
\BD
  \sum_{\btop{j=1}{j\:odd}}^{q-1}\chi(k)^j
\ED
vanishes if $\chi^2(k)\ne 1$. Indeed, if we put $\zeta=\chi(k)\ne\pm 1$, we have
\BD
  \sum_{\btop{j=1}{j\:odd}}^{q-1}\zeta^j=\sum_{j=0}^{q-1}\zeta^j-\sum_{j=0}^{q/2-1}\zeta^{2j}.
\ED
Each of the sums on the right hand side equals $0$. Hence (\ref{2.8}) follows, and (\ref{2.6}) gives the assertion of Theorem \ref{t1}.

\MN
{\em Remark.} By (\ref{2.4}),
\BD
  \sum_{\btop{k=1}{\chi^2(k)=1}}^{p-1}k\chi(k)\equiv p\cdot l\equiv 1 \MB{ mod }2.
\ED
In particular, $E\ne 0$.

\bigskip
\centerline{\bf Competing interests and data availability}

\MN
The author declares that there are no competing interests. The paper has no associated data.


\MN
Institut f\"ur Mathematik \\
Universit\"at Innsbruck   \\
Technikerstr. 13/7        \\
A-6020 Innsbruck, Austria \\
Kurt.Girstmair@uibk.ac.at

\end{document}